\documentclass{amsart}
\usepackage{geometry}
\usepackage{amsthm,amsmath,amssymb,amsxtra,latexsym,mathrsfs,url}

\usepackage{ab}

\usepackage{cite} 

\usepackage{graphicx} %

\theoremstyle{plain}
\newtheorem{proposition}{Proposition}[section] 
\newtheorem{lemma}[proposition]{Lemma} 
\newtheorem{theorem}[proposition]{Theorem} 
\newtheorem{corollary}[proposition]{Corollary}
\newtheorem{method}[proposition]{Method}
\newtheorem{conjecture}[proposition]{Conjecture}
\newtheorem{proposition-short}{Proposition} 
\newtheorem{theorem-short}{Theorem} 

\theoremstyle{remark}
\newtheorem{definition}{Definition} 
\newtheorem{example}{Example}
\newtheorem{examples}{Examples}
\newtheorem{remark}{Remark}

\newcommand{\thm}{\begin{theorem}}
\newcommand{\Thm}{\end{theorem}}
\newcommand{\thma}{\begin{theorem*}}
\newcommand{\Thma}{\end{theorem*}}

\newcommand{\conj}{\begin{conjecture}}
\newcommand{\Conj}{\end{conjecture}}
\newcommand{\prf}{\begin{proof}}
\newcommand{\Prf}{\end{proof}}
\newcommand{\dfn}{\begin{definition}}
\newcommand{\Dfn}{\end{definition}}
\newcommand{\prop}{\begin{proposition}}
\newcommand{\Prop}{\end{proposition}}
\newcommand{\cor}{\begin{corollary}}
\newcommand{\Cor}{\end{corollary}}
\newcommand{\cora}{\begin{corollary*}} 
\newcommand{\Cora}{\end{corollary*}}
\newcommand{\lem}{\begin{lemma}}
\newcommand{\Lem}{\end{lemma}}
\newcommand{\lema}{\begin{lemma*}}
\newcommand{\Lema}{\end{lemma*}}
\newcommand{\rmk}{\begin{remark}}
\newcommand{\Rmk}{\end{remark}}
\newcommand{\exa}{\begin{example}}
\newcommand{\Exa}{\end{example}}
\newcommand{\exas}{\begin{examples}}
\newcommand{\Exas}{\end{examples}}
\newcommand{\meth}{\begin{method}}
\newcommand{\Meth}{\end{method}}
\newcommand{\bexa}{\begin{barredexample}}
\newcommand{\bExa}{\end{barredexample}}

\newcommand{\set}[1]{\{#1\}}
\newcommand{\lies}{\subset} 
\newcommand{\eand}{\text{ and }} 
\newcommand{\eimplies}{\text{ implies }} 
\newcommand{\id}{\text{id}} 
\newcommand{\of}{\circ} 
\newcommand{\call}[1]{\mathop{\forall #1 \,}} 
\DeclareMathOperator{\dom}{dom} 
\DeclareMathOperator{\cod}{cod} 
\newcommand{\sous}[1]{{}_{#1}} 
\newcommand{\iso}{\cong} 
\newcommand{\Z}{\mathbb{Z}}
\newcommand{\Set}{\ensuremath{\textrm{\upshape{{\bf Set}}}}} 
\DeclareMathOperator{\Ob}{Ob} 
\DeclareMathOperator{\Mor}{Mor} 
\DeclareMathOperator{\Nat}{Nat} 
\newcommand{\natto}{\Rightarrow} 
\DeclareMathOperator{\colim}{colim} 
\newcommand{\dleq}{\sqsubseteq}

\newcommand{\bottom}{\bot} %
\newcommand{\vertof}{\mathbin{\raisebox{0.2ex}{$\scriptscriptstyle{\Delta}$}}} 
\DeclareMathOperator{\troot}{root} 
\DeclareMathOperator{\tright}{right} 
\DeclareMathOperator{\tleft}{left} 
\DeclareMathOperator{\height}{height} 
\DeclareMathOperator{\source}{{\bf s}} 
\DeclareMathOperator{\target}{{\bf t}} 

\newenvironment{rightbrace}
  {\left.\begin{aligned}}
  {\end{aligned}\right\rbrace}

\usepackage{scalerel,stackengine}
\stackMath
\newcommand\verywidehat[1]{%
\savestack{\tmpbox}{\stretchto{%
  \scaleto{%
    \scalerel*[\widthof{\ensuremath{#1}}]{\kern-.6pt\bigwedge\kern-.6pt}%
    {\rule[-\textheight/2]{1ex}{\textheight}}%
  }{\textheight}%
}{0.5ex}}%
\stackon[1pt]{#1}{\tmpbox}%
}

\usepackage{amsxtra} 
\stackMath %
\newcommand\verywidetilde[1]{ 
\savestack{\tmpbox}{\stretchto{ %
  \scaleto{ %
    \scalerel*[\widthof{\ensuremath{#1}}]{\kern-0.1pt \sptilde \kern-0.3pt} %
    {\rule[-\textheight/2]{1ex}{\textheight}} %
  }{\textheight} %
}{1.5ex}} %
\stackon[-2pt]{#1}{\tmpbox} %
}

\begin{document}
\title{Foundations of a Recent Extension of Category~Theory and Topos~Theory}
\author{
Lucius T. Schoenbaum
}
\date{\today}
\keywords{category, topos, categorical logic, programming language theory, type theory}


\begin{abstract}
This article is an introduction to the basic generalized category theory used in the recent work \cite{ScMd} studying an extension 
of the theory of categories and categorical logic, including parts of topos theory. %
We discuss functors, equivalences, natural transformations, adjoints, and limits in a generalized setting, %
giving a concise outline of these frequently arising constructions. 
\end{abstract}

\maketitle

\section{Introduction}\label{s.i}

Category theory \cite{MacCW,KaSc1,BaWe1} has its origins in mathematics, and has since been applied to analysis of mathematical foundations and programming languages. It begins with the insight that diagrams and morphisms have mathematical properties that are independent of function theory, and independent of any use of points as arguments. %
In  \cite{EiMa1}, the paper on natural transformations in which the elementary notions of category theory are introduced for the first time, Eilenberg and MacLane write:
\begin{quote}
It is thus clear that the objects play a secondary role, and could be entirely omitted from the definition of a category. However, the manipulation of the applications would be slightly less convenient were this done. 
\end{quote}
Therefore we can say that there are two views which have been known to category theorists since the very beginning of the subject: %
a one-sorted definition describing a universe of pure maps, and a two-sorted definition including the objects that, in applications, are prior to the maps that they inspire. %
These two views pull against one another in a way that seems perhaps like a natural, irresolvable tension. %
The latter approach has proven to be the dominant one, while the former approach has made occasional appearances, for example in work by Ehresmann \cite{Ehresmann1}, Street \cite{StT1}, and %
in recent work by %
Cockett \cite{CockettConstellations}. %

The two-sorted view masks the potential for generalization that begins with the (less often used) one-sorted formulation. In the latter case, an axiom requires that the source and target maps $s$ and $t$ are trivial upon iteration: $ss = st = s, tt = ts = t$. This condition, however, is extraneous, as it never arises in proofs. Dropping it gives rise to a rather general notion, %
which may be %
weakened further via replacing some equalities with inequalities, as suggested by some kinds of applications \cite{SmPl1,Reynolds1980}. %
In \cite{ScMd} the theory of categories and categorical logic is developed in this generalized setting, including parts of topos theory. %
Interested readers can consult \cite{ScMd} as a more extended reference, and find there some discussion of applications. %
This short article is an introduction to the basic generalized category theory used in \cite{ScMd}: the theory of functors, equivalences, natural transformations, adjoints, and limits in the generalized setting. %

%

%

%

\section{Generalized Categories}\label{s.gencat}

{\em Preliminaries.} 
We write composition $G \of F := (f \mapsto G(F(f)))$ and in general, for mappings $F$ and $G$ with common domain and codomain (i.e., in which concatenation is meaningful) we define the operation
$$G \vertof F := (f \mapsto G(f)F(f)),$$
the standard vertical composition operation \cite{MacCW}. %
In any context where it is meaningful, we use the standard arrow notation $F:A \to B$ to mean that an element $F$ is given, the domain of $F$ is $A$, and the codomain of $F$ is $B$. %
The notation $\downarrow$ (cf. Definition \ref{d.gencat})
indicates that all composed pairs of elements in the expression or relation are in fact composable pairs. 
A given mapping $F$ equipped with a domain $\dom(F)$ and a codomain $\cod(F)$, which (especially when regarded as formal constructs) may be denoted using the following conventions:
$$\dom(F) = \source(F) = \bar{F},$$
$$\cod(F) = \target(F) = \hat{F}.$$

%
\subsection{The 2-Category of Categories}\label{ss.2catofcat}

Before starting on the generalization of category theory that is the focus of this article, we offer a brief review of perhaps the most important elementary construction in category theory: the strict 2-category of categories. %
This section assumes some familiarity with category theory \cite{MacCW}. %

Let $\sC, \sD$ be categories. Two natural transformations $\beta: G \natto H, \alpha: F \natto G$ between functors $F,G: \sC \to \sD$ may be composed via the rule
$$\beta \vertof \alpha (X) := \beta(X) \cdot \alpha(X)$$
where $(\cdot)$ denotes composition in $\sD$. This gives a category $\Nat(\sC, \sD)$. Identities in $\Nat(\sC, \sD)$ are given by $\id_F (X) := \id_X$. 

Given natural transformations $\alpha: F \natto G$ between functors $\sC \to \sD$, and $\beta: F' \natto G'$ between functors $\sD \to \sE$, we obtain a well-defined function $\Ob(\sC) \to \Mor(\sE)$ via
$$\beta \star \alpha (X) := \alpha ( \hat \beta (X) ) \cdot \bar \alpha ( \beta (X)),$$
where hats and bars are used as defined in section \ref{s.gencat} below. This can also be written
$$\beta \star \alpha = (\alpha \of \hat\beta) \vertof (\bar\alpha \of \beta)$$
Note that 
$$\bar \alpha (X) = \overline{\alpha(X)},$$
$$\hat \alpha (X) = \widehat{\alpha(X)}.$$

\prop\label{p.fivefacts}
In the notation above, whenever expressions on both sides of the formula are defined, we have:
\begin{enumerate}
	\item $\beta \star \alpha = (\hat\alpha \of \beta) \vertof (\alpha \of \bar\beta).$ %
	\item $\beta \star \alpha$ is a natural transformation $G \of F \natto G' \of F'.$ 
	\item $(\gamma \star \beta) \star \alpha = \gamma \star (\beta \star \alpha).$
	\item If
		$$\begin{rightbrace}
			\alpha: F \natto G \\
			\beta: G \natto H 
		\end{rightbrace} : \sC \to \sD, $$
		$$\begin{rightbrace}
			\alpha' : F' \natto G' \\
			\beta' : G' \natto H' 
		\end{rightbrace} : \sD \to \sE, $$
		then
		$$(\beta' \vertof \alpha') \star (\beta \vertof \alpha) = (\beta' \star \beta) \vertof (\alpha' \star \alpha).$$
	\item If $\id_F^{\vertof}$ is the identity of $F$ with respect to the product $\vertof$ in $\Nat(\sC, \sD)$, then 
		$$\alpha \star \id^{\vertof}_F = \alpha,$$
		$$\id^{\vertof}_F \star \beta = \beta,$$
		whenever both sides are defined.
\end{enumerate}
\Prop
\prf
(1) 
\begin{align*}
	(\beta \star \alpha) (X)		&= (\alpha \of \hat \beta) \vertof (\bar \alpha \of \beta) (X) \\
						&= \alpha(\hat \beta(X)) \cdot \bar \alpha (\beta(X)) \\
						&= \alpha(\widehat{\beta(X)}) \cdot F(\beta(X)) \\
						&= G(\beta(X)) \cdot \alpha(\overline{\beta(X)}) \\
						&= \hat \alpha (\beta(X)) \cdot \alpha(\bar \beta (X)) \\
						&= (\hat \alpha \of \beta) \vertof (\alpha \of \bar \beta) (X).
\end{align*}

(2) by Fact 1.

(3) Apply the definition.

(4) by Fact 1 and since $\hat \alpha = \bar \beta, \widehat{\alpha'} = \overline{\beta'}$.

(5) direct calculation.
\Prf

{\em Remarks.} 
\begin{enumerate}
	\item We may write simply $\id_F$ or $1_F$ in light of Fact (5).
	\item Fact 4 is often referred to as the {\em interchange law}.
\end{enumerate}

An immediate consequence of Proposition \ref{p.fivefacts} is the following: {\em The category of categories is a strict two-category.} %
By ``the category of categories'' is meant the set of small categories, functors, and natural transformations in a fixed universe $\sU_{univ}$.

\subsection{Definition}\label{ss.gencat}

\dfn\label{d.gencat}
A {\em generalized category} is a structure $(\sC, \dleq, \source,\target, \cdot)$ where $\sC$ is a set, $\dleq$ is a relation on $\sC$, $\source$ and $\target$ are mappings $\sC \to \sC$, and $(\cdot)$ is a partially defined mapping $\sC \times \sC \to \sC$, denoted $a \cdot b$ or $ab$. These are required to satisfy
\begin{enumerate}
	\item $(\sC, \dleq)$ is a partially ordered set, 
			\label{ax.gencat-po}
	\item $ab$ $\downarrow$ if and only if $\source(a) \dleq \target(b)$. \label{ax.gencat-po-comp}
	\item If $(ab)c$ $\downarrow$ or $a(bc)$ $\downarrow$ then $(ab)c = a(bc)$. \label{ax.gencat-assoc}
	\item If $ab$ $\downarrow$ then $\source(ab) = \source(b)$ and $\target(ab) = \target(a)$. \label{ax.gencat-comp-st}
	\item (Element-Identity) For all $a \in \sC$, there exists $b \in \sC$ such that \label{ax.gencat-element-id}
		\begin{enumerate}
			\item $\source(b) = \target(b) = a$,
			\item if $bc$ $\downarrow$ then $bc = c$,
			\item if $cb$ $\downarrow$ then $cb = c$,
		\end{enumerate}
	\item (Object-Identity) Let $a \in \sC$ and $\source(a) = \target(a) = a$. Then \label{ax.gencat-object-id}
		\begin{enumerate}
			\item if $ba$ $\downarrow$ then $ba = b$. 
			\item If $ab$ $\downarrow$ then $ab = b$.
		\end{enumerate}
	\item (Order Congruences\footnote{These axioms are needed for the Kleisli construction \cite{ScMd}.})
		\begin{enumerate}
			\item If $a \dleq b$ then $\source(a) \dleq \source(b)$ and $\target(a) \dleq \target(b)$. \label{ax.gencat-order1}
			\item $a \dleq b$ and $c \dleq d$ and $ac,bd$ $\downarrow$ implies $ac \dleq bd.$ \label{ax.gencat-order2}
			\item $a \dleq b$ implies $1_a \dleq 1_b$. \label{ax.order3}
		\end{enumerate}
\end{enumerate}
The element $c$ of axiom (\ref{ax.gencat-element-id}) %
is unique, and is denoted $1_a$ or $\id_a$, and called the {\em identity} on $a$. %
\Dfn

As a partially ordered set a generalized category resembles, but is weaker than, a domain \cite{GrZ+}, indeed motivation for the ordering comes from domain theory \cite{SmPl1, WaD1}. If $a \dleq b$, we say that $a$ {\em approximates} $b$ or {\em upcasts to} $b$, and $b$ {\em sharpens} $a$ or {\em downcasts to} $a$. %
When the ordering $\dleq$ is nontrivial, one may call $\sC$ a {\em casting} generalized category. %
We often think of casting categories as having at least a bottom element $\bottom$, but we do not assume this in the definition, since we would like, as a special case, for an ordinary one-category to be a generalized category. %
If the order given by $\dleq$ is discrete, we might say that the generalized category is {\em discrete}, and similarly for other order-theoretic attributes, but as this may lead to confusion with the notion of a discrete category (one with essentially no morphisms), we shall say instead that such a generalized category is a {\em sharp} generalized category. We allow ourselves to refer to a {\em casting} generalized category whenever we wish to emphasize that we refer to a generalized category that is not assumed to be sharp. 

An {\em element} $f \in \sC$ is an element $f$ of the underlying set $\sC$. 
An {\em object} $a$ in $\sC$ is an element $a$ of $\sC$ such that $\source(a) = \target(a) = a$. 
We write $\Ob(\sC)$ for the set of objects. %
For $a \in \sC$, we define the {\em height} of $a$, denoted $\height(a),$
to be the maximum of the set of nonnegative integers $n$ such that there exists a sequence $\vec s$ of source and target operations of length $n$ such that $\vec s (i)$ is an object, %
unless there is an infinite sequence $\vec s$ of source and target operations such that no subsequence yields an object. In that case, we say that $\height(a) = \infty$. 

With this terminology, Definition \ref{d.gencat} says that in a generalized category with identities, every element $a$ has an identity $1_a$, and that if the element is an object, this identity is $a$ itself. %
If $a \in \sC$ has identity $1_a$ and is not an object, then $a \neq 1_a$. %

The maps $\source$ and $\target$ of the definition are called the  {\em source} or {\em domain} and {\em target} or {\em codomain} maps, respectively. As noted above, we allow ourselves for convenience to denote the map $\source(a)$ by either $\dom(a)$ or $\bar a$, and the map $\target(a)$ by either $\cod(a)$ or $\hat a$. %

Given a generalized category $\sC$, any element of $\sC$ may be composed with other compatible elements, and it is equipped with a ``tail'' of fellow elements, defined by the $\source$ and $\target$ maps. We think of the product %
as developing from right to left, and we may write $c:a \to b$
when $\source(a) = b$, $\target(a) = c$. %
Note as an aside that if one pictures instead a representation $a = \sous{c}a_b$ of $a$, one has a picture of composition $\sous{c}a_b \,\sous{b}d_e = \sous{c}(ad)_e$. %
This notation can be iterated to 
$$a = \sous{\sous{g}c_f}a_{\sous{e}b_d}$$ 
In this manner one can visualize a binary tree. %

\subsection{An Alternative Approach}\label{ss.secondapproach}

We now proceed to define a generalized category using an alternative approach, and discuss why we choose the approach of Definition \ref{d.gencat}. %

\dfn\label{d.gencatalt}
A {\em generalized category} is a structure $(\sC, \dleq, \source,\target, \cdot)$ where $\sC$ is a set, $\dleq$ is a relation on $\sC$, $\source$ and $\target$ are operators (mappings) $\sC \to \sC$, and $(\cdot)$ is a partially defined binary operation $\sC \times \sC \to \sC$, denoted $a \cdot b$ or $ab$. These are required to satisfy
\begin{enumerate}
	\item $(\sC, \dleq)$ is a partially ordered set, %
	\item If $(ab)c$ $\downarrow$ or $a(bc)$ $\downarrow$ then $(ab)c = a(bc)$. 
	\item If $ab$ $\downarrow$ then $\source(ab) = \source(b)$ and $\target(ab) = \target(a)$.
	\item $ab$ $\downarrow$ if and only if $\source(a) \dleq \target(b)$.
	\item (Object-Identity) Let $a \in \sC$ and $\source(a) = \target(a) = a$. Then \label{ax.gencat-objid}
		\begin{enumerate}
			\item if $ba$ $\downarrow$ then $ba = b$. 
			\item If $ab$ $\downarrow$ then $ab = b$.
		\end{enumerate}
	\item (Order Congruences)
		\begin{enumerate}
			\item If $a \dleq b$ then $\source(a) \dleq \source(b)$ and $\target(a) \dleq \target(b)$. \label{ax.gencat-order1}
			\item $a \dleq b$ and $c \dleq d$ and $ac,bd$ $\downarrow$ implies $ac \dleq bd.$ \label{ax.gencat-order2}
			\item $a \dleq b$ implies $1_a \dleq 1_b$. \label{ax.gencat-order3}
		\end{enumerate}
\end{enumerate}
A generalized category is said to be equipped {\em with identities} if for every $a \in \sC$, if there exists $b \in \sC$ such that $\source(b) = a$ or $\target(b) = a$, then there exists $c \in \sC$ such that $c b$ $\downarrow$ implies $c b = b$, and $b c$ $\downarrow$ implies $b c = b$.%
The element $c$ %
is unique, and is denoted $1_a$ or $\id_a$, and called the {\em identity} on $a$. %
An {\em element} $f \in \sC$ is an element $f$ of the underlying set $\sC$. %
An {\em object} $a$ in $\sC$ is an element $a$ of $\sC$ such that $\source(a) = \target(a) = a$. %
A {\em subject} $U$ in $\sC$ is an element $U$ of $\sC$ such that there exists $f \in \sC$ such that $\source(f) = U$, or there exists $f \in \sC$ such that $\target(f) = U$. %
\Dfn

The approach of Definition \ref{d.gencat} has the advantage of having fewer basic concepts than Definition \ref{d.gencatalt}. %
All elements are subjects and all elements have identities. %
This makes many steps of the development go smoothly. %
On the other hand, Definition \ref{d.gencat} creates so many identities that one sometimes wonders if they are better avoided after all. %
Thus one might seem to be at an impasse concerning whether Definition \ref{d.gencat} or Definition \ref{d.gencatalt} is more preferable. %
This ambivalence is resolved by the notion of ideal category \cite{ScMd}. %
Ideal categories arise naturally in categorical logic, %
and they may be computationally implemented. %
In such categories, and in particular in the generalized category of contexts $\boC\Lambda$, there are identities present %
just as Definition \ref{d.gencat} requires. 

\subsection{Resuming, from Definition \ref{d.gencat}}

\prop\label{p.duality}
Up to reversal of $\dleq$, Definition \ref{d.gencat} is symmetric in the source and target maps $\source$ and $\target$. Therefore every proof $\Phi$ about a generalized category $\sC$ continues to hold when, in all assumptions, definitions, and deduction steps, composition, the order $\dleq$, and the role of source and target are reversed. 
\Prop

Such a proof $\Phi'$ is said to be obtained from $\Phi$ ``by duality'' \cite{MacCW}. %
This simple fact has a profound effect on the entire subject. %
The generalized category formed by the operation of Proposition \ref{p.duality} is called the {\em opposite generalized category} $\sC^{op}$ of $\sC$. %

\begin{example}
Let $\sC$ be a category \cite{MacCW}. %
Then {\em the generalized category generated by $\sC$} is obtained from $\sC$ by identifying the identity $1_X$ of each object $X \in \sC$ with $X$, and closing over $1_{()}$. %
Considering a concrete example, such as the generalized category generated by the category of all groups, we may write $\id_X$ for $X$, with the identification $\id_X = X$ being understood. %
More formally, we define: %

\dfn\label{d.cat}
A generalized category $\sC$ is a {\em category} or {\em one-category} if the source and target of every nonidentity $f$ in $\sC$ is an object in $\sC$. 
\Dfn

We now have a rough ontology:

$$
\begin{tabular}{c|c}
\begin{tabular}{c}
sharp category \\
= category 
\end{tabular} 
	& 
\begin{tabular}{c}
casting category
\end{tabular} 
			\\ \hline
\begin{tabular}{c}
sharp generalized category 
\end{tabular} 
	& 
\begin{tabular}{c}
casting generalized category \\
= generalized category \\
\end{tabular}
\end{tabular}
$$
\end{example}

\begin{example}
In some instances it is possible to write down a generalized category explicitly. There is an empty generalized category, and $\sC = \set{a:a \to a}$, the trivial generalized category. %
More generally, any set $S$ is a generalized category after setting $\source(a) = \target(a) = a$ for $a \in S$, we say that the generalized category is {\em discrete} or a {\em zero-category}, or simply that it is a set. %
(Thus, sets and categories are examples of generalized categories.) 
Because of the identity axiom, other than finite sets there are no finite generalized categories. 
To amend language, we therefore define:

\dfn\label{d.figen}
A generalized category $\sC$ is {\em finitely generated} if there is a finite set $\sC$ such that the remainder of $\sC$ consists only of identities. 
\Dfn

There are many examples of generalized categories that are not ordinary categories, the simplest perhaps being
$\sC = \set{a: a \to a, b: a \to b}$.
Another simple example is 
$\sC = \set{a: b \to b, b: a \to a}$.
This generalized category is finite, but does not possess objects, moreover every element is a subject. %
A generalized category may also lack objects due to infinite descent, for example
$\sC = \set{a_n : a_{n-1} \to a_{n-1} \mid n \in \Z}.$ %
\end{example}

\begin{example}
Motivation for the casting relation $\dleq$ in a casting generalized category comes from %
the subtyping relation in some type theoretical systems \cite{Reynolds1980,Pierce1}. Subtyping is a feature found in many programming languages, including %
most (if not all) object-oriented languages, which typically involves some form of field/method inheritance. %
Another commonplace form of subtyping is the explicit and implicit typecasting of built in types, for example to treat a single-precision integer as a double precision one. %
Implementing subtyping involves data type {\em coercion,} or modification of a type at compile time or at run time. %
Type-theoretically, condition (\ref{ax.gencat-comp-st}) of Definition \ref{d.gencat} corresponds to a type system for a language which allows automatic upcasting upon evaluation at subtypes. %
An algorithm for typecasting that is sensitive to the input could be implemented using a variant of dependent typing mechanisms.
\end{example}

\begin{example}\label{ex.binarytree}
Let $\sC$ be a generalized category, and consider 
the condition on $\sC$ that hom sets should contain a unique element or else be empty. 
To obtain a (possibly infinite) planar binary tree one adds the condition that source and target may not loop except trivially, that is, for every element $a \in \sC$, and for every finite sequence $(x_1, \dots, x_n)$ where $x_i$ is either $\source$ or $\target$ (source or target) if $x_n x_{n-1} \dots x_1 a = a$ then it is required that $\source a = \target a = a$, that is, or (using the terminology of trees) that $a$ is a leaf.  %
Presheaves on such trees arise for example in database theory, see for example \cite{SpK3}. %
\end{example}

\exa
A {\em generalized (directed) graph} \cite{ScMd} is simply a triple $(\sA, \source, \target)$, where $\sA$ is a carrier set, %
and $\source, \target$ are maps $\sA \to \sA$. %
An element of $\sA$ is (synonymously) an {\em edge}. %
An {\em object} in a generalized graph is an element $a \in \sA$ such that $\source a = \target a = a$, %
that is, a common fixed point of the endomorphisms $\source$ and $\target$. 
Ordinary graphs correspond bijectively with 1-dimensional generalized graphs, 
where we say that generalized graph is {\em 1-dimensional} if $\source \source = \source \text{ and } \target \target = \target.$
With the obvious composition via compound paths, a generalized graph becomes a (sharp) generalized category. 

There are plentiful settings where generalized graphs may arise. For example, suppose that there is a system of goods $\sA_0$. %
The edges of $\sA$ are certificates (issued, say perhaps, by different governing bodies) that say that a good $a \in \sA_0$ may be exchanged for another good $b \in \sA_0$. %
Suppose it is accepted that a good is always exchangeable for itself. %
Now let's suppose that such certificates themselves may be exchanged, %
but that this requires that one has a higher-level certificate for this higher-level trade. %
If we imagine a certain impetus exists among those we imagine making the exchanges, we can expect that there will next arise trading for these certificates as well, 
giving rise to a generalized graph (in fact, a generalized deductive system, via a simple extension of Kolmogorov's reasoning about intuitionistic logic in \cite{Kolmogorov1}). 
\Exa

\exa\label{x.setf}
For a planar binary tree $\ft$, let
\begin{align*}
\troot(\ft) &\text{ is the root of $\ft$.} \\
\tleft(\ft) &\text{ is the tree given by the left descendant of the root, and its descendants.} \\
\tright(\ft) &\text{ is the tree given by the right descendant of the root, and its descendants.}
\end{align*}
From any category $\sC$ we can form a sharp generalized category $\sC f$ as follows: take the set $\sC f$ to be the set of all planar binary trees of morphisms in $\sC$, subject to a source-and-target condition %
$$\dom \troot(\dom \ff) = \dom \troot(\ff) ,$$
and
$$\cod \troot(\cod \ff) = \cod \troot(\ff),$$
where if $\ff$ be such a tree,
$$\cod \ff = \tleft(\ff),$$
the left descendent tree of $\ff$, and 
$$\dom \ff = \tright(\ff),$$
the right descendent tree of $\ff$. 
These conditions set up a recursive condition on elements of $\sC f$.
For $\fg, \ff \in \sC f$, we set
$$\fg \cdot \ff := (\text{ the tree $\fh$ with left descendent $\troot(\ff)$, right descendent $\troot(\fg)$ and root $\troot(g)\cdot \troot(f)$. }).$$
This is a well-defined product, by the source-and-target condition above. It is checked that this is a (sharp) generalized category. An element of $\sC f$ may be visualized as
$$
\raisebox{-0.5\height}{\includegraphics{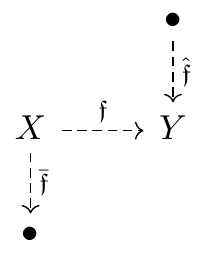}}
$$
Constructions on the original $\sC$ can be carried over to $\sC f$, for example, if $\sC$ has products (equalizers, coproducts, coequalizers), then so (respectively) does $\sC f$. If $\sC$ is (co)complete, however, it does not imply that $\sC f$ is (co)complete, see \cite{ScMd}. 

Lawvere's comma category construction \cite{LwV63, MacCW} may also be observed to yield generalized categories, even when the input data is an ordinary category. %
Fix two generalized categories $\sC$, $\sD$, and $\sE$, and functors $S:\sD \to \sC$ and $T: \sE \to \sC$. %
Let 
$$
(S,T) := \set{(d,e,\ff) \mid d \in \sD, e \in \sE, \ff \text{ is a planar binary tree of pairs $(f,g)$, $f,g \in \sC$ }}
$$
Set
$$\overline{(d,e,\ff)} = (\bar d, \bar e, \tright(\ff)),$$
$$\widehat{(d,e,\ff)} = (\hat d, \hat e, \tleft(\ff)).$$
Composition in $(S,T)$ is defined as in $\sC f$.
\Exa

%

%

%

\section{Elementary Theory, Category of Invertibles}\label{s.elem}

We now define functors and hom sets: 

\dfn\label{d.functor}
A mapping $\sC \to \sC'$ between generalized categories is {\em functorial} or a {\em functor} if 
\begin{enumerate}
	\item $a \dleq b \eimplies F(a) \dleq F(b),$
	\item $F(\bar a) = \overline{F(a)},$
	\item $F(\hat a) = \widehat{F(a)},$
	\item $F(ab) = F(a)F(b), \text{ if } ab \,\downarrow,$
	\item $F(1_a) = 1_{F(a)}.$
\end{enumerate}
We thus have a category ${\bf GenCat}$ of generalized categories and functors. %
\Dfn

Functors are also called {\em covariant functors}. A {\em contravariant functor} from $\sC$ to $\sC'$ is a unital map satisfying
\begin{enumerate}
	\item if $a \dleq b$ then $F(b) \dleq F(a)$,
	\item $F(\bar a) = \widehat{F(a)},$
	\item $F(\hat a) = \overline{F(a)},$
	\item $F(ab) = F(b)F(a) \text{ if } ab \,\downarrow,$
\end{enumerate}
instead of the corresponding covariant relations. %

\dfn\label{d.homsets}
The sets 
$$\hom(a,b) = \set{c \in \sC \mid \bar c = a, \hat c = b},$$ 
for $a,b \in \sC$, are called the {\em hom sets} of $\sC$. 
\Dfn

\dfn
A {\em subcategory} of a generalized category $\sC$ is a subset $\sC'$ of $\sC$ whose order is inherited from $\sC$ closed under source, target, composition, and identities: if $a \in \sC'$, then $1_a \in \sC'$. %
A subcategory $\sC'$ is {\em full} if $a,b \in \sC'$ implies $\hom(a,b)$ is contained in $\sC'$. 
\Dfn

The composition of two functors is a functor, and functors send objects to objects. 

\dfn
Two generalized categories $\sC$ and $\sC'$ are {\em isomorphic} if there is an invertible functor (i.e., invertible as a mapping) $F$ from $\sC$ to $\sC'$. 
\Dfn

\prop\label{p.flattening}
There is a functor, flattening, from the category of generalized categories to the category of categories. 
\Prop
\prf
Let $\sC$ be a generalized category with identities. %
Let $\Ob(\sC_{flat})$ be $\set{[f] \mid f \in \sC}$, the objects of $\sC$ indexed by the elements of $\sC$. 
Let $\Mor(\sC_F)$ again be a set $\set{(f) \mid f \in \sC}$ indexed by the elements of $\sC$, and define source and target
$$s((f)) = [s(f)],$$
$$t((f)) = [t(f)].$$
Then $\sC_{flat}$ is a category whose composition and identities are
$$(g) \cdot (f) := (gf),$$
$$
1_{[f]} = (1_f).
$$
Given a functor $F: \sC \to \sD$ in ${\bf GenCat}$, we immediately obtain a functor $\sC_{flat} \to \sD_{flat}$. 
\Prf

Note that $\sC_{flat}$ contains a flattening of the identity structure, even in cases where $\hom(a,a) = \set{1_a}$.  

There is also a category $flat \sC$, the further flattening of $\sC$ to a zero-category. %
It is defined by:
$$flat(f) := \begin{cases} 	(f), 	& \text{if $f = 1_g$ for some $g \in \sC$, } \\ 
					[f]	& \text{otherwise,}
	\end{cases}
$$
where $[f]$ is defined by $\source([f]) = \target([f]) = [f]$, and $(f): flat(\source(f)) \to flat(\target(f)).$

\dfn\label{d.invertible}
If $\sC$ is a generalized category, an element $a \in \sC$ is {\em invertible} if there exists $b \in \sC$ such that $ab = 1_{\hat a}$ and $ba = 1_{\bar a}$. %
\Dfn

\prop\label{p.invertible}
$\phantom{v}$
\begin{enumerate}
	\item The inverse $a^{-1}$ of an element $a$ of $\sC$ is unique if it exists. 
	\item $\widehat{a^{-1}} = \bar a$ and $\overline{a^{-1}} = \hat a$. (Even if $\sC$ is casting.) 
	\item All objects $a$ are invertible: $a^{-1} = a$. 
	\item Functors send invertibles to invertibles: $F(\theta^{-1}) = F(\theta)^{-1}$.
\end{enumerate}
\Prop

There are a few ways a generalized category may be partitioned into equivalence classes:

\dfn\label{d.epimonicisoclass}
For $a, b \in \sC,$ we have the following equivalence relations:
\begin{enumerate}
	\item $a$ and $b$ are in the same {\em monic class}, or {\em subobject}, $a \sim_m b$, if there exists invertible element $\theta \in \sC$ such that $a\theta = b$. 
	\item $a$ and $b$ are in the same {\em epic class}, or {\em quotient}, $a \sim_e b$, if there exists invertible element $\theta \in \sC$ such that $\theta a = b$; 
	\item $a$ and $b$ are in the same {\em iso class}, $a \sim b$, if there exist invertible elements $\theta_1, \theta_2 \in \sC$ such that $\theta_1 a = b \theta_2$. 
\end{enumerate}
\Dfn

Let $\Theta$ denote the set of all invertible elements in $\sC$. Define the symbol 
$$a\Theta := \set{a\cdot\theta \mid \theta \in \Theta \eand a\cdot\theta \downarrow},$$
and define the symbols $\Theta a, \Theta a\Theta$, etc. similarly. Then for $a,b \in \sC$,
$b$ belongs to the monic class of $a$ if and only if $b \in a\Theta$, %
$b$ belongs to the epic class of $a$ if and only if $b \in \Theta a$, %
and 
$b$ belongs to the iso class of $a$ if and only if $b \in \Theta a \Theta$. %
This notation is useful for back-of-the-envelope calculations, but it can be misleading: it need not be true that $\Theta f \Theta = \Theta g \Theta$, even if $f$ and $g$ are invertible. %

\dfn\label{d.epimoniciso}
An element $m$ of a generalized category $\sC$ is {\em monic} if $mf,mg \downarrow$ and $mf = mg$ implies $f = g$. %
An element $e$ in $\sC$ is {\em epi} if $fe,ge \downarrow$ and $fe = ge$ implies $f = g$. %
We say $a$ is {\em isomorphic} to $b$, denoted 
$$a \iso b,$$
if there exists an invertible element $\theta$ with $\bar\theta = a, \hat\theta = b$. %
\Dfn

If $a$ is monic and $a \sim_m b$, then $b$ is monic, and the $\theta$ given by the definition is unique. Similarly, if $a$ is epic and $a \sim_e b$. 

For every $a,b \in \sC$, $a$ is isomorphic to $b$ iff $1_a$ is in the same iso class as $1_b$, that is,
$$a \iso b \quad \Longleftrightarrow \quad 1_a \sim 1_b.$$
For $a,b$ objects, this becomes: 
$$a \iso b \quad \Longleftrightarrow \quad a \sim b.$$

\prop\label{p.isoskel}
Let $\sC$ be a generalized category. Then the set of iso classes forms a sharp category. %
The objects of this category are the iso classes of invertible elements of $\sC$. 
\Prop
\prf
Let $\tilde \sC$ be the set of iso classes of $\sC$, %
let $\tilde a, \tilde b$, ... denote elements in $\tilde \sC$. 
Define
$$\tilde a \cdot \tilde b := %
\set{\theta_1 a \theta_2 b \theta_3 \mid \theta_1,\theta_2,\theta_3 \text{ invertible, and $\theta_1 a \theta_2 b \theta_3 \, \downarrow$}}.$$
This is a partially defined map $\tilde \sC \times \tilde \sC \to \tilde \sC$. For $a \in \sC$, let 
$$\tilde \source (\tilde a) := \widetilde{1_{\source a}},$$
$$\tilde \target (\tilde a) := \widetilde{1_{\target a}}.$$
These operations are well-defined: if $a = \theta_1 b \theta_2$, then $\bar a$ is isomorphic to $\bar b$, so, say, $\theta 1_{\bar b} \theta^{-1} = 1_{\bar a}$, so $\widetilde{1_{\source a}} = \widetilde{1_{\source b}}$, and similarly for $\tilde \target$. %

We take the order $\dleq$ on $\tilde \sC$ to be trivial, and we check Definition \ref{d.gencat}. %
The first four conditions are immediate: for (\ref{ax.gencat-comp-st}),  %
if $\tilde a, \tilde b \in \tilde \sC$, then $ \tilde a \tilde b \,\downarrow$. This occurs 
if and only if $ \set{\theta \in \sC \mid \theta: \bar a \to \hat b \text{ is invertible}}$ is nonempty, 
if and only if $ 1_{\bar a} \sim 1_{\hat b}$, 
if and only if $\tilde \source(\tilde a) = \tilde \target(\tilde b)$. %
Next, we observe that if $\tilde a$ is an element of the form $\tilde\source \tilde b$ or $\tilde \target \tilde b$ 
in $\tilde \sC$, then it must be of the form $\widetilde{1_b}$ for some $b \in \sC$, and 
$$\tilde\target(\widetilde{1_b}) = \tilde \source(\widetilde{1_b}) = \widetilde{1_b},$$
so $\tilde{1_b}$ is an object. Next, we have
$$\widetilde{1_a} \cdot \tilde b = \set{\theta_1 1_a \theta_2 b \theta_3 } = \set{\theta_4 b \theta_3} = \tilde b,$$
and similarly, $\tilde b \widetilde{1_a} = \tilde b$ whenever the product is defined. %
So $\tilde \sC$ is a sharp generalized category, in fact a one-category, after closing over $1_{()}$. %
The second statement is merely the observation that $a$ is invertible if and only if $\tilde a = \widetilde{1_{\source a}} = \widetilde{1_{\target a}}.$
\Prf

\dfn\label{d.catinvertibles}
We refer to the category $\tilde \sC$ of Proposition \ref{p.isoskel} as the {\em category of invertibles} of $\sC$. %
\Dfn

The {\em skeleton} of a generalized category $\sC$ is any full subcategory such that each element of $\sC$ is isomorphic in $\sC$ to exactly one element of the subcategory. %
Skeletons are unique up to isomorphism \cite{MacCW}. %
In the case of a category $\sC$, %
the category of invertibles expresses exactly the same data as a skeleton, but in a different way: any iso class that is an object in the category of invertibles contains not a set of invertibles in $\sC$ that are pairwise isomorphic, but instead, the set of all the isomorphisms that relate them pairwise to one another. %
On the other hand, an iso class that is an arrow in the category of invertibles is a noninvertible arrow $f \in \sC$ well-defined up to a commutative square with invertible columns. 

Since every element has an identity, thus taking the category of invertibles is the same as the operation of flattening (Proposition \ref{p.flattening})
followed by taking the skeleton, yielding the description just made in the previous paragraph. Thus it is perhaps natural to think of it as the ``category of identities'' of the generalized category. 

It is also the case that a functor $F$ lifts to a functorial map $\tilde F$ on the category of invertibles.
Indeed, define 
$$\tilde F : \tilde \sC \to \tilde \sC',$$
via 
$$\tilde F(\tilde a) := \widetilde{F(a)}.$$
This is well-defined, as a consequence of (2) (which depends on the unital property of $F$):
$$\tilde F (\theta_1 a \theta_2) = \verywidetilde{F(\theta_1) F(a) F(\theta_2)} = \widetilde{F(a)}.$$
So we check functoriality: we have
$$\widetilde{1_{\source(F(a))}} = \widetilde{1_{F(\source(a))}} = \widetilde{F(1_{\source(a)})} = \tilde{F}(\widetilde{1_{\source(a)}}) = \tilde{F}(\source(\tilde a)),$$
and
$$\widetilde{1_{\source(F(a))}} = \source(\widetilde{F(a)}) = \source(\tilde{F}(\tilde{a})).$$
Similarly, 
$$\target(\tilde{F}(\tilde{a})) = \tilde{F}(\target(\tilde{a})).$$
And
$$\tilde{F}(\tilde{a}\tilde{b}) = \tilde{F}(\widetilde{a \theta b}) = \widetilde{F(a\theta b)} = \verywidetilde{F(a) \theta' F(b)} = \tilde{F}(\tilde{a}) \tilde{F}(\tilde{b}).$$
Finally, $\tilde{F}$ is unital since $\tilde \sC$ and $\tilde{\sC'}$ are categories.

A notion weaker than isomorphism arises from considering the categories of invertibles. 

\dfn\label{d.equivalent}
Generalized categories $\sC$ and $\sC'$ are {\em equivalent} if their categories of invertibles are isomorphic. 
\Dfn

This definition appeals directly to a comparison of the categories of invertibles. 
Now consider two functors $F,G: \sC \to \sC'$ that both define the same functor $\tilde{\sC} \to \widetilde{\sC'}$ on the categories of invertibles of $\sC$ and $\sC'$. This can only mean that there exist a pair of functions $\theta_1,\theta_2:\sC \to \sC'$ such that $\call{a \in \sC} \theta_i(a)$ is invertible for $i=1,2$, and for all $a \in \sC$, %
$$ \theta_1 (a) F(a) = G(a) \theta_2 (a) \,\downarrow.$$
If this holds we may write 
$$F \iso G.$$

\prop\label{p.equivalent}
Two generalized categories $\sC$ and $\sC'$ are equivalent if either of the following two equivalent conditions are satisfied. 
\begin{enumerate}
	\item Their categories of invertibles are isomorphic via a pair $\tilde F, \tilde G$, where $\tilde G = \tilde F^{-1}$, that come from functors
	 $F: \sC \to \sC'$ and $G: \sC' \to \sC$. 
	\item There exist two functors $F, G$ from $\sC \to \sC'$ ($\sC' \to \sC$, respectively) satisfying 
		$$F \of G \iso \id_{\sC'},$$
		$$G \of F \iso \id_{\sC}.$$
\end{enumerate}
\Prop

We can consider properties that a functor $\tilde F$ on the category of invertibles has as an ordinary functor, and view them as properties of the underlying functor $F$:

\dfn\label{d.essinjsurj}
A functor $F:\sC \to \sC'$ is {\em essentially injective} if it satisfies one of the following equivalent conditions, 
\begin{enumerate}
	\item $\tilde F$ is injective. 
	\item For $a,b \in \sC$, $F(a) = F(b)$ implies $a \sim b$.
\end{enumerate}
and $F$ is {\em essentially surjective} if it satisfies one of the following equivalent conditions: 
\begin{enumerate}
	\item $\tilde F$ is surjective.
	\item For $\alpha \in \sC'$, there exists $a \in \sC$ with $F(a) \sim \alpha$.
\end{enumerate}
\Dfn

From our initial investigation of equivalences between generalized categories, we arrived at the notion of equivalence via a pair of functors $F$ and $G$. We could, however, view this machinery (the pair ($\theta_1, \theta_2$)) as instead relating the two functors, and extend it:

\dfn\label{d.mof}
Let $\sC, \sC'$ be generalized categories, let $F,G:\sC \to \sC'$ be two functors. We say that a {\em morphism of functors} \cite{KaSc1} from $F$ to $G$ is a pair $(\theta_1,\theta_2)$ of maps $\sC \to \sC'$ satisfying, for all $a \in \sC$,
\begin{equation}\label{eq.mof}
\theta_1 (a) F(a) = G(a) \theta_2 (a) \, \downarrow
\end{equation}
Note that here, $\theta_1$ and $\theta_2$ are no longer presumed to be invertible. 
We may write the morphism of functors with the notation 
$(\theta_1, \theta_2): F \Rightarrow G$. 
\Dfn

Note that the maps $\theta_1$ and $\theta_2$ are maps from $\sC$ to $\sC'$, not from $\Ob(\sC)$ to $\sC'$ (cf. \cite{MacCW}). %

\begin{example}
Let $A = (a_{ij})$ be a matrix with coefficients in a ring $R$, and let $f:R \to S$ be a ring homomorphism. One naturally sets $f(A) = (f(a_{ij}))$, and doing this, one sees that
\begin{equation}\label{eq.fdet}
\det(f(A)) = f(\det(A)).
\end{equation}
This relation can be interpreted by observing that $GL_n$ is a functor from the category of rings to the category of groups, and likewise for the mapping that sends a ring to its group of units, and a ring homomorphism to the pointwise-identical homomorphism on the respective groups of units. So if $f:R \to S$, and writing $F(f)$ for the map defined above extending $f$ to a map on $GL_n (R)$, and $G(f)$ for the map changing $f$ to a map on the group of units, we have
$$\det() \of F(f) = G(f) \of \det()$$
by rewriting equation (\ref{eq.fdet}). From this expression we can read off the morphism of functors: 
$$\theta_1(f) = \det : GL_n(S) \to S^{\times},$$
$$\theta_2(f) = \det : GL_n(R) \to R^{\times}.$$
We see that in this example, $\theta_1$ and $\theta_2$ come from a single map $\theta$ on the objects (rings). This is not only typical of categories, it is guaranteed to happen. Indeed, if we return to the general situation of Definition \ref{d.mof}, inserting $a = 1_b$ into equation (\ref{eq.mof}) gives 
$$ \theta_1 (1_b) = \theta_2 (1_b)$$
for $b \in \sC$, so in particular, for all objects $b$,
$$ \theta_1 (b) = \theta_2 (b).$$
Thus $\theta_1$ and $\theta_2$ are identical on objects, and since one-categories have no higher morphisms, this single map on objects completely characterizes $(\theta_1, \theta_2)$. %
\end{example}

In the terminology of section \ref{s.nat} that follows, this means that a morphism of functors between functors relating categories is always {\em natural}. %
In the setting of generalized categories, we might suppose that this naturality property is a condition special to one-categories, %
since it does not appear to have any a priori motivation. %
However, %
the theory that results from dropping the naturality condition appears to be significantly weaker: %
\begin{enumerate}
	\item There is no strict 2-category of non-natural transformations, functors, and generalized categories. 
		Here, the wheel turns on the tiniest of pedestals: in the notation of \ref{ss.2catofcat}, the relations
		$$\bar \alpha (X) = \overline{\alpha(X)},$$
		$$\hat \alpha (X) = \widehat{\alpha(X)}$$
		hold only in the natural setting. So we do not prove (1) of Proposition \ref{p.fivefacts}. 
	\item While there is a notion of non-natural adjunction, there is no hom set bijection. %
		A key step in the proof uses the naturality of the unit and counit maps. %
		This in turn is used to prove that left adjoints are right exact. %
	\item Because there is no adjoint hom set bijection, some theorems relating equivalences of categories with 
		properties of functors %
		no longer hold. In particular a full, faithful, essentially surjective functor might not define an equivalence. 
\end{enumerate}
For these reasons, we do not take the development any further until we introduce naturality in the next section.

%

%

%

\subsection{Globular Sets}\label{ss.globularsets}

This section is about the relationship between generalized categories and globular sets \cite{BaN1}. 
A globular set is a presheaf of shape $\mathbb{G}$ (that is, a functor $\mathbb{G}^{\text{op}} \to \Set$), where $\mathbb{G}$ is the category of natural numbers $n \geq 0$ together with maps
$$
\raisebox{-0.5\height}{\includegraphics{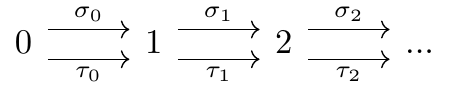}}
$$
subject to the relations 
$
\sigma_{i+1} \of \sigma_{i} = \tau_{i+1} \sigma_{i}, \quad \tau_{i+1} \of \tau_{i} = \sigma_{i+1} \of \tau_{i},
$
for $i \geq 0$.

\dfn
Let $\sC$ be a generalized category. A {\em $k$-cell} in $\sC$ is an element $f$ of $\sC$ such that for every $k$-element sequence $\vec s$ of operations $\source$ and $\target$ that satisfy when applied to $f$, %
\begin{enumerate}
	\item $\source^k f$ and $\target^k$ are objects, and $\source^{k-n} f$ and $\target^{k-n}$ are not objects, for all $0 \leq n \leq k$, 
	\item $\source \target f = \source \source f$ and $\target \source f = \target \target f$,
	\item $\source f$ and $\target f$ is are $k-1$-cells.
\end{enumerate}
For example, in a 1-dimensional category, all elements are 1-cells, and some elements are also 0-cells. 
An element $f$ of a generalized category $\sC$ is {\em cellular} if $f$ is a $k$-cell for some $k \geq 1$, %
and a generalized category $\sC$ is {\em cellular} if every element of $\sC$ is cellular. 
\Dfn

\prop
There is an equivalence (given by a forgetful-free adjunction) between sharp, cellular generalized categories and the category of globular sets. 
\Prop
\prf
To prove this, we must be sure clarify the statement: when referring to sharp, cellular generalized categories, we refer not to the full subcategory but to the category whose morphisms $F: \sC \to \sD$ are subject to the extra condition 
\begin{enumerate}
	\item for all $a \in \sC$, $\source (F(a)) = F(a)$ implies $\source a = a$.
\end{enumerate}
This says we cannot map $k$-cells for $k > 0$ to $0$-cells. %
Then let $\dim(a) := \min\set{n \mid \source^n a = \source^{n+1} a }$. %
Define a mapping 
$$\sC \mapsto (n \mapsto \set{a \in \sC \mid \dim a = n}).$$ %
to the category of globular sets, for a sharp cellular generalized category $\sC$. This is the desired equivalence. 
\Prf

Examples of noncellular generalized categories are abundant, for example arising from the theory of trees and related notions, see for example \cite{ElBlTi1}.

%

%

%

\section{Naturality}\label{s.nat}

In this section we establish the second of the two notions of equivalence we consider, namely
natural equivalence. %
As already noted, {\em the distinction between natural and non-natural vanishes in the case of categories.} %
Under natural equivalence, we obtain a 2-category of generalized categories, and in particular, an interchange law (Theorem \ref{t.functorcategory}). We can also establish, using the final lynchpin that naturality provides so to speak, the hom set bijection associated with adjoint pairs (Theorem \ref{t.nathombij}). Consequently the familiar rule that an equivalence between categories is given by a fully faithful essentially surjective functor carries over to generalized categories (Theorem \ref{t.ffes}). The full and faithful properties are tied to the naturality condition, which gives rise to maps not only on individual elements, but on entire hom sets. 

\dfn\label{d.nattrans}
Let $\sC, \sD$ be generalized categories, let $F,G: \sC \to \sD$. %
Let $(\theta_1, \theta_2) :F \Rightarrow G$ be a morphism of functors. %
We say that $(\theta_1, \theta_2)$ is {\em natural} or that $(\theta_1, \theta_2)$ is a {\em natural transformation} if, 
for every $a,b \in \sC$, %
$$\theta_1(a) = \theta_1 (b)$$
whenever $\hat a = \hat b$, and
$$\theta_2 (a) = \theta_2 (b)$$ 
whenever $\bar a = \bar b$.
\Dfn

Thus, naturality means that the function $\theta_1(a)$ can be replaced with the function $\hat a \mapsto \theta_1(1_{\hat a})$ of the element $\hat a$, and $\theta_2$ can be replaced with the function $\bar a \mapsto \theta_2(1_{\bar a})$ of the element $\bar a$. But, as noted in section \ref{s.elem}, 
$\theta_1(1_b) = \theta_2(1_b)$ 
for all elements $b$. 
Hence a natural transformation reduces to a single map $\theta: \sC \to \sC'$, from which $\theta_1$ and $\theta_2$ are immediately derived: 
$$\theta_1 (a) := \theta(1_{\hat a}),$$
$$\theta_2 (a) := \theta(1_{\bar a}).$$
We refer to a natural transformation $(\theta_1, \theta_2)$ by referring to this map $\theta$. In terms of $\theta$ the defining relation of a morphism of functors becomes
$$\theta(\hat f \hspace{0.8pt}) \cdot F(f) = G(f) \cdot \theta(\bar f \hspace{0.8pt}) \, \downarrow.$$

\dfn\label{d.naturaleq}
Two generalized categories $\sC$ and $\sC'$ are {\em naturally equivalent} %
 if they are equivalent via natural transformations
$$\theta:F\of G \iso \id_{\sC'},$$
$$\theta':G\of F \iso \id_\sC.$$
\Dfn

Naturally equivalent generalized categories are, in particular, equivalent (Definition \ref{d.equivalent}). %
With the extra condition of naturality, 
the way is clear to extend many 
justly well-known results of one-category theory \cite{MacCW} to the generalized setting: 

\thm\label{t.functorcategory}
The system given by all of the 
generalized categories, 
functors, 
and natural tranformations
forms a strict 2-category. %
\Thm
\prf
We define the products 
$$\theta_1 \vertof \theta_2,$$
$$\theta_1 \star \theta_2$$ 
just as in section \ref{ss.2catofcat}, and proceed as in the one-categorical case. %
\Prf

We include the naturality condition when defining adjoints: %

\dfn\label{d.adjunction}
Let $\sC$ and $\sD$ be generalized categories. An %
{\em adjunction} %
$(F,G,\eta,\varepsilon)$ is a pair of functors 
\vspace{-1ex}
$$
\raisebox{-0.5\height}{\includegraphics{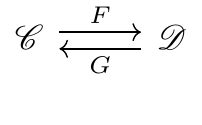}}
$$

\vspace{-15pt}
\noindent together with natural transformations
$$
\eta: \id_\sC \to G\of F, \quad \varepsilon: F \of G \to \id_\sD,
$$
satisfying the identities
\begin{align}\label{eq.triid}
(G \of \varepsilon) \vertof (\eta \of G) &= 1_G, \\
(\varepsilon \of F) \vertof (F \of \eta) &= 1_F,
\end{align}
where $1_F$ is the mapping $f \mapsto 1_{F(f)}$. %
Given an adjunction $(F,G,\eta,\varepsilon)$, $\eta$ is called the {\em unit} and $\varepsilon$ is called the {\em counit} of the adjunction. 
A natural equivalence $(\theta, \theta')$ is an {\em adjoint equivalence} if $\theta$ and $\theta'$ are the unit and counit of an adjunction. 
\Dfn

\thm\label{t.nathombij}
Let $\sC, \sD$ be generalized categories, and let $F,G: \sC \to \sD$ be functors. The following are equivalent:
\begin{enumerate}
	\item $(F,G,\eta,\varepsilon)$ forms an adjunction \!\!\!
$
\raisebox{-0.6\height}{\includegraphics{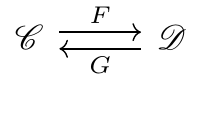}}
$ \!\!\!.

\vspace{-12pt}
	\item %
For every $f$ in $\sC$ and $g$ in $\sD$, there is a %
bijection of sets
\begin{equation}\label{eq.adjbij}
\hom(F(f), g) \iso \hom(f, G(g)),
\end{equation}
that is natural in $f$ and $g$. This means that if $\phi_{f,g}$ is the bijection (\ref{eq.adjbij}), then 
for every $k: g \to g'$, and $h:f' \to f$, the following diagrams commute: 
$$
\raisebox{-0.5\height}{\includegraphics{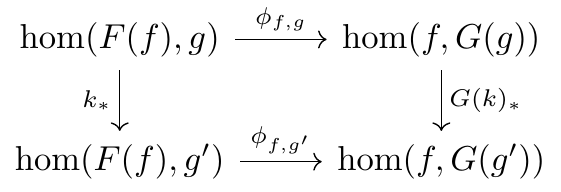}}
\quad\quad
\raisebox{-0.5\height}{\includegraphics{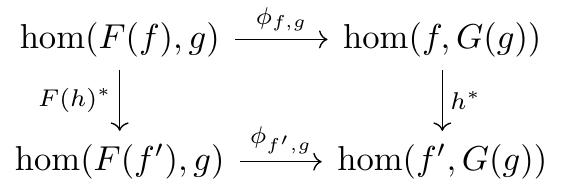}}
$$
Equivalently $\phi$ satisfies
\[
u \cdot F(v): F(f) \to g \eimplies \phi(u \cdot F(v)) = \phi(u) \cdot v,
\]
\[
v' \cdot v: F(f) \to g \eimplies \phi(v' \cdot v) = G(v') \cdot \phi(v).
\]
\end{enumerate}
\Thm
\prf
The proof is formally the same as in the one-categorical case (see \cite{MacCW}). 
\Prf

\dfn\label{d.faithfulfull}
Let $\sC,\sD$ be generalized categories, $F: \sC \to \sD$ a functor. 
For $a,b \in \sC$, let $F_{a,b}$ be the mapping on the domain $\hom(a,b)$ given by $f \mapsto F(f)$. %
We say that $F$ is {\em faithful} if %
for all $a,b$, $F_{a,b}$ is injective, 
and 
we say that $F$ is {\em full} if %
for all $a,b$, $F_{a,b}$ is surjective. 
\Dfn

Thus for example full means: if $\alpha, \beta$ in $\sD$ are of the form $F(a), F(b)$, for $a,b \in \sC$, and if $\gamma: \alpha \to \beta$, then $\gamma$ is of the form $F(c)$ for $c \in \sC$. 

\thm\label{t.ffes}
Let $\sC, \sD$ be generalized categories, and let $F: \sC \to \sD$ be a functor. The following are equivalent:
\begin{enumerate}
	\item $F$ is a natural equivalence,
	\item $F$ is a natural adjoint equivalence,
	\item $F$ is full, faithful, and essentially surjective.
\end{enumerate}
\Thm
\prf
The proof, much the same as in the one-categorical case, is left to the reader. %
\Prf

%

%

%

\section{Limits}\label{s.lim}

In this section we establish the elements of the theory of limits and colimits in sharp generalized categories. 
We consider limits with respect to mappings $I \to \sC$ as in Definition \ref{d.limit} that are weaker than functors. This, for example, allows us to form the shape of a product or coproduct of any set of elements in a generalized category. %

\dfn\label{d.limitfunctor}
Let $\sC, \sC'$ be generalized categories. A {\em functor up to objects} from $\sC$ to $\sC'$ is a map $F:\sC \to \sC'$ satisfying, for every $a,b \in \sC$,
\begin{enumerate}
	\item $F(ab) = F(a)(b)$,
	\item $F(a)$ is an identity in $\sC'$ if and only if $a$ is an identity in $\sC$,
	\item $F(\source(a)) = \source(F(a))$ unless $a$ is an object of $\sC$,
	\item $F(\target(a)) = \target(F(a))$ unless $a$ is an object of $\sC$. 
\end{enumerate}
\Dfn

\dfn\label{d.cone}
Let $\sC$ be a generalized category, $I$ a generalized category (the index of a cone needs only be a set, but in practice it is always a (generalized) category). A {\em cone} in $\sC$ with index $I$ is a map $\sigma:I \to \sC$ such that 
$$\text{for all $i,j \in I$, } \overline{\sigma(i)} = \overline{\sigma(j)}.$$ %
Dually, {\em cocone} in $\sC$ with index $I$ is a map $\sigma:I \to \sC$ such that %
for all $i,j \in I$, $\widehat{\sigma(i)} = \widehat{\sigma(j)}.$ %
A cone or cocone is {\em finitely generated} if the index set $I$ is finitely generated (Definition \ref{d.figen}). %
This common source is the {\em vertex} of the cone, and the {\em vertex} of a cocone is the common target. %
Given a cone or cocone $\pi$, we may refer to $\pi(i)$ for some $i \in I$ as a {\em member} of the cone. %
\Dfn

\dfn\label{d.limit}
Let $\sC, I$ be generalized categories. Let $\alpha: I \to \sC$ be a functor, possibly only a functor up to objects. %
A cone is said to be {\em over (or below) the base $\alpha$} if %
\begin{enumerate}
	\item $\widehat{\pi(i)} = \alpha(i)$, for all $i \in I$, %
	\item for all $i \in I$, $\pi(\hat i) = \alpha(i) \pi(\bar i)$. %
\end{enumerate}
A {\em limit} of $\alpha$ is a cone $\pi:I \to \sC$ below the base $\alpha$ such that for any cone $\tilde \pi:I \to \sC$ over the same base $\alpha$, there is a unique $\lambda \in \sC$ such that $\tilde \pi = \pi \vertof \lambda$. (Here, $\pi \vertof \lambda$ is the map defined by $(\pi \vertof \lambda)(i) = \pi(i) \cdot \lambda.$) %

Dually, a cocone is said to be {\em over (or below) the base $\alpha$} if 
\begin{enumerate}
	\item $\overline{\pi(i)} = \alpha(i)$, for all $i \in I$,
	\item for all $i \in I$, $\pi(\bar i) = \pi(\hat i) \alpha(i)$.
\end{enumerate}
A {\em colimit} of $\alpha$ is a cocone $\pi:I \to \sC$ such that for any cone $\tilde \pi:I \to \sC$ over the base $\alpha$, there is a unique $\lambda \in \sC$ such that $\tilde \pi = \lambda \vertof \pi$. Here, $\lambda \vertof \pi$ is the map defined by $(\lambda \vertof \pi)(i) = \lambda \cdot \pi(i)$, as before. 
\Dfn

Thus a cone fits a pattern as in the following Figure:
$$
\raisebox{-0.5\height}{\includegraphics{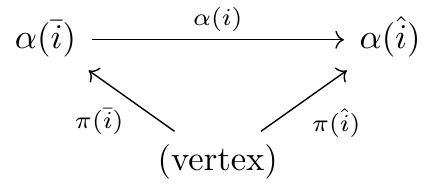}}
$$

\noindent 
The word limit is often used to refer to the domain of the cone, and similarly colimit is used to refer to the codomain of the cocone. %
The terms {\em product}, {\em equalizer}, {\em coproduct}, {\em coequalizer}, etc. retain their meaning from ordinary categories, referring to limits based on diagrams $\alpha: I \to \sC$ of the same shape as in the one-categorical case, and where $\alpha$ may be a functor only up to objects. %
We follow standard terminology and say that a generalized category {\em has finite limits} if there is a limit cone for every finitely generated diagram $\alpha: I \to \sC$, and dually for colimits. 

We denote the set of limits of the functor $\alpha:I \to \sC$ by $\lim(\alpha,I)$ or just $\lim \alpha$. We denote the colimit $\colim(\alpha,I)$ or simply $\colim(\alpha)$. %

If $\sC$ is a generalized category, there exist (finitely generated) diagrams $J \to \sC$ that cannot be defined and do not exist in an ordinary category. 
However, we still have:

\thm\label{t.prodeqsuffice}
Let $\sC$ be a generalized category. For $\sC$ to have all finite limits, it suffices that $\sC$ has all finite products and equalizers. 
\Thm
\prf
We proceed by induction on the height of finitely generated diagrams $\alpha: I \to \sC$. 
A finitely generated diagram of height $0$ is a finite product, hence it has a limit cone in $\sC$ by hypothesis. %
Suppose that all finitely generated diagrams of height $k \geq 0$ have a limit cone, and 
let $\alpha: I \to \sC$ be a diagram of height $k + 1$. %
Define 
$$\alpha^{\leq k}$$
to be $\alpha$ restricted to the generalized category $I^{\leq k}$ formed by taking the collection of all elements of $I$ of height $\leq k$, along with all identities of $I$. 
It is easy to see that $I^{\leq k}$ is closed under composition, thus it is a generalized category. Therefore $\alpha^{\leq k}$ is a diagram on $\sC$, and by hypothesis, has a limit cone $\sigma^{\leq k}$ with vertex, say, $L^{\leq k}$. %
Consider $flat(I^{\leq k})$, the flattening of $I^{\leq k}$ to a zero-category (section \ref{s.elem}). %
The diagram $flat(\alpha^{\leq k}): flat(I^{\leq k}) \to \sC$ induced by $\alpha^{\leq k}$ %
is a diagram of height zero, so it has a limit cone $\sigma^{\leq k, flat}$, with vertex, say, $L^{\leq k, flat}$. 
The cone $\sigma^{\leq k}$ on $I^{\leq k}$ induces a cone on $flat(I^{\leq k})$, 
so there exists a universal arrow 
$$u_1: L^{\leq k} \to L^{\leq k, flat}.$$
Now let $I^{k+1, flat}$ be the flattened (to a zero category) elements of $I$ of height $k+1$. %
The diagram $\alpha$ induces a diagram $\alpha^{k+1, flat}$ on $I^{k+1, flat}$, %
defined by 
$$\alpha^{k+1, flat} (i) := \target(\alpha(i)).$$
This diagram (of height zero) has a limit cone $\sigma^{k+1, flat}$ 
with vertex, say, $L^{k+1, flat}$. %
For $i \in I$ of height $k+1$, let $\pi_i$ be the element in $\sC$ which is the projection 
$$\pi_i : L^{\leq k, flat} \to \target(\alpha(i)),$$
coming from the diagram $\sigma^{\leq k, flat}$ on $I^{\leq k, flat}$ (where our notation hides this fact about $\pi_i$). 

The previous cone $\sigma^{\leq k, flat}$ with vertex $L^{\leq k , flat}$ itself has projection arrows to the elements $\target(\alpha(i))$ 
as $i$ ranges over $\alpha^{k+1, flat}$. 
Therefore, there is a universal arrow 
$$u_2: L^{\leq k, flat} \to L^{k+1, flat}.$$
Moreover, for each $i$ of height $k+1$, there is also a projection arrow to the element $\source(\alpha(i))$, and 
composing each of these projection arrows with $\alpha(i)$ gives a second cone with the same vertex $L^{\leq k, flat}$ %
on the diagram $\alpha^{k+1, flat}$. 
So we may again find a universal arrow
$$u_3: L^{\leq k, flat} \to L^{k+1, flat},$$
by applying the universal property of the limit with vertex $L^{k+1, flat}$ a second time. %
We compose $u_2$ and $u_3$ with $u_1$ to form parallel arrows, and take the equalizer:
\[
\raisebox{-0.5\height}{\includegraphics{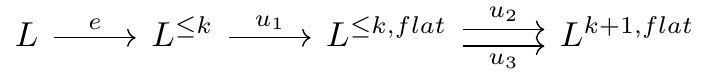}}
\]
Now we define, for $i$ in $I$ of height $\leq k+1$,
$$\sigma^{\leq k+1} (i) := \pi_i \cdot u_1 \cdot e.$$
We claim that this is a limit cone for the diagram $\alpha^{\leq k+1}: I^{\leq k+1} \to \sC$. 
Since we pass through $e$ to reach $L^{\leq k+1}$, $\sigma^{\leq k+1}$ satisfies 
$\sigma^{\leq k+1} (\hat i) = \alpha^{\leq k+1} (i) \cdot \sigma^{\leq k+1} (\hat i)$, 
hence is a limit cone. 
Suppose that 
$\tilde \sigma^{\leq k+1}: I^{\leq k+1} \to \sC$ is a diagram with vertex, say, $\tilde L$ satisfying 
$\tilde \sigma^{\leq k+1} (\hat i) = \alpha^{\leq k+1} (i) \tilde \sigma^{\leq k+1} (\bar i).$ 
Then $\tilde \sigma^{\leq k+1}$ restricts to a cone on $\alpha^{\leq k}$, hence there is a universal arrow
$$\tilde e: \tilde L \to L^{\leq k}.$$
Because $\tilde \sigma^{\leq k+1}$ has the limit property even at the height $k+1$, $\tilde \sigma^{\leq k+1}$ satisfies $u_2 \cdot u_1 \cdot \tilde e = u_3 \cdot u_1 \cdot \tilde e$, and thus $\tilde e$ factors through $e$ uniquely, as desired. 
\Prf

\dfn
Let $F: C \to C'$ be a functor. Then $F$ {\em preserves limits} or is {\em left exact} if for every functor $\alpha:I \to C$, 
$$F(\lim(\alpha)) \lies \lim(F\of \alpha).$$
Dually, $F$ {\em preserves colimits} or is {\em right exact} if for every functor $\alpha:I \to C$,
$$F(\colim(\alpha)) \lies \colim(F\of \alpha).$$
$F$ is said to {\em create limits} if for every element $\pi \in \lim(F \of \alpha)$, there exists a unique $\pi' \in \lim(\alpha)$ such that $F(\pi') = \pi$. Dually, $F$ is said to {\em create colimits} if for every element $\pi \in \colim(F \of \alpha)$, there exists a unique $\pi' \in \colim(\alpha)$ such that $F(\pi') = \pi$. %
\Dfn

For example, the hom functor 
$$b \mapsto \hom(-,b)$$
preserves limits. %
Dually, the contravariant hom functor 
$$a \mapsto \hom(a,-)$$
preserves colimits. These functors may be extended to generalized categories \cite{ScMd}. 

\thm\label{t.hlalex}
Let $F: \sC \to \sD$ be a functor between generalized categories $\sC$ and $\sD$. 
Then if $F$ has a left adjoint $G: \sD \to \sC$, then it is left exact. 
\Thm
\prf
Like the proof for categories, the proof for generalized categories relies on naturality of the adjoints via the bijection (\ref{eq.adjbij}). 
\Prf

The dual statement to \ref{t.hlalex} is immediate: a functor with a right adjoint is right exact. %

%

%

%

\section{Conclusion}\label{s.conc}

We have surveyed the beginning of the formal investigation in \cite{ScMd} of assumptions about the basic notions in category theory. 
There are numerous advanced notions of category theory that have not yet made an appearance in our development, 
for example, ends, coends, monads, Kan extensions, double categories, toposes (elementary or sheaf-theoretic), to name a few. %
Some of these are treated in \cite{ScMd}, so the interested reader can consult there for details. 

We have observed that there exists a theory of functors, natural transformations, adjoint pairs, limits, and colimits for generalized categories. %
More precisely, there are (cf. section \ref{s.gencat}) two generalizations that are combined into one larger one:  %
First, by allowing an approximate operation of composition, 
and second, by allowing generalized higher cells. 
We have seen that the structure of limits and natural transformations is similar to the structure as it arises in ordinary one-categories, 
so that the casting structure does not stymie the theory's formal development. %
This suggests that a subtyping structure may appear on the semantic side, as it often does on the type theoretic side %
when there is a programming language of interest. %
We have investigated a notion of non-natural transformation suggested by the one-categorical case where naturality is not a necessary assumption, and we have found that the device of non-natural equivalence is not sufficiently strong. Therefore, we have argued that naturality must be an explicit assumption in the generalized setting. 

\bibliographystyle{abbrv}
\bibliography{mathmain}

\vspace*{10pt}
\noindent {\footnotesize LUCIUS T. SCHOENBAUM \\
\hspace*{18pt}\address{DEPARTMENT OF MATHEMATICS AND STATISTICS\\
\hspace*{18pt}UNIVERSITY OF SOUTH ALABAMA\\
\hspace*{18pt}MOBILE, AL 36688-0002}\\
{\it E-mail}: schoenbaum@southalabama.edu\\
}
\clearpage

\end{document}